\documentclass[12pt]{amsart}
\usepackage{amssymb,amsmath,a4wide}

\newtheorem{theorem}{Theorem}
\newtheorem{lemma}[theorem]{Lemma}
\newtheorem{corollary}[theorem]{Corollary}

\newtheorem{remark}{\it Remark\/}

\newcommand{\pa}{\partial}
\newcommand{\ep}{\varepsilon}
\newcommand{\RR}{{\mathbb{R}}}
\renewcommand{\div}{\mathop{{\rm div}}}
\newcommand{\supp}{\mathop{{\rm supp}}}

\begin{document}

\title[Renormalized variational principles]{Renormalized variational principles and Hardy-type inequalities}%
\author{Satyanad Kichenassamy}%
\address{Laboratoire de Math\'ematiques (UMR 6056), CNRS \&\ Universit\'e
de Reims Champagne-Ardenne, Moulin de la Housse, B.P. 1039,
F-51687
Reims Cedex 2\\France}%
\email{satyanad.kichenassamy@univ-reims.fr}%

\thanks{Appeared in {\bf Differential and Integral Equations}, 19 : 4 (2006) 437-448}%

\date{August 2005}%
\begin{abstract}
Let $\Omega\subset\RR^2$ be a bounded domain on which Hardy's
inequality holds. We prove that $[\exp(u^2)-1]/\delta^2\in
L^1(\Omega)$ if $u\in H^1_0(\Omega)$, where $\delta$ denotes
the distance to $\pa\Omega$. The corresponding higher-dimensional
result is also given.
These results contain both Hardy's and Trudinger's
inequalities, and yield a new variational
characterization of the maximal solution of the Liouville equation
on smooth domains, in terms of a renormalized functional. A global
$H^1$ bound on the difference between the maximal solution and the
first term of its asymptotic expansion follows.
\end{abstract}
\maketitle

\section{Introduction}
\label{sec:intro}

\subsection{Hardy's and Trudinger's inequalities}

Let $\Omega$ be an arbitrary domain in $\RR^N$, $N\geq 2$; let
$\delta(x)$ denote the distance of $x$ from $\pa\Omega$.

If $\Omega$ is bounded with Lipschitz boundary, and $u\in
H^1_0(\Omega)$, the generalized Hardy's inequality states
that
\[ \| \frac{u}{\delta}\|_{L^2(\Omega)}\leq H\|\nabla
u\|_{L^2(\Omega)}.
\]
The optimal value of the ``Hardy constant'' $H$, as well as
possible generalizations and improvements of this inequality, have
been the object of much attention, see \cite{MMP,Mz,OK}. Hardy's
inequality arises naturally in several variational problems of
recent interest, as well as in the proof of decay estimates
\cite{BM,SK-div}.

On the other hand, if $N=2$, Trudinger's inequality
implies that $\exp(u^2)-1$ is integrable. This
suggests our first result:
\begin{theorem}\label{th:1-2d}
Let $\Omega\subset\RR^2$ be such that Hardy's
inequality holds. Then, for any $u\in H^1_0(\Omega)$,
\[[\exp(u^2)-1]/\delta^2\in L^1(\Omega).
\]
\end{theorem}
This theorem will be derived from a result of independent
interest:\footnote{The Lebesgue measure $dx$ is understood in all
integrals in this paper.}
\begin{theorem}\label{lem:1-2d}
Let $\Omega\subset\RR^2$ be an arbitrary domain
with $\pa\Omega\neq\emptyset$.
Then, if $u\in H^1_0(\Omega)$, and $q>2$,
\[\left(\int_\Omega\frac{u^q}{\delta^2}\right)^{1/q} \leq \Sigma_q
\left(\int_\Omega |\nabla u|^2+\frac{u^2}{\delta^2}\right)^{1/2},
\]
where $\Sigma_q=O(q^{1/2+1/q})$ as $q\to\infty$.
\end{theorem}
\begin{remark}
These results admit natural generalizations to higher dimensions,
which are stated and proved
in section \ref{sec:2}.
Th.~\ref{lem:1-2d} may be considered as trivially true if
$\pa\Omega=\emptyset$ with the convention that
$\delta\equiv+\infty$ in this case.
For background results on Trudinger's
inequality, see \cite{T,M-70,A,BW}; we merely note that our
argument is closer to Trudinger's than to Moser's, because the
distance function does not transform in a convenient manner under
symmetrization.
\end{remark}
\begin{remark}
No regularity or boundedness assumptions on $\Omega$ are required.
This is somewhat surprising in view of the fact \cite[p.~120]{A}
that elements of $H^1(\Omega)$ are not necessarily $L^q$ for
$q>2$, if $\Omega$ is unbounded and with finite volume. Note that
for domains with thin ``ends'' at infinity, $\delta$ is very
small, and the r.h.s.~of our inequality is not equivalent to the
$H^1$ norm.
\end{remark}
In the situation of Th.~\ref{th:1-2d}, we find in particular that
\begin{equation}\label{eq:exp}
\frac{e^{2u}-1-2u}{\delta^2}\in L^1(\Omega).
\end{equation}
As an application of this result, we present next a variational
characterization of the maximal solution of an elliptic problem
with monotone nonlinearity. The question of constructing
a renormalized energy in this context was raised by H. Brezis last
year in the Nonlinear Analysis Seminar of Paris VI.

\subsection{A renormalized energy for boundary blow-up}

Let $\Omega\subset\RR^2$ be a bounded domain of class
$C^{2+\alpha}$, with $0<\alpha<1$. The distance function
$\delta(x)$ is of class $C^{2+\alpha}$ near and up to the
boundary, but is only Lipschitz over $\Omega$ in general. It is
therefore convenient to introduce a function $d(x)\in
C^{2+\alpha}(\overline\Omega)$, which coincides with $\delta(x)$
near $\partial\Omega$, and is positive inside
$\Omega$.\footnote{One may simply take $d=F(\delta)$ for an
appropriate $F$.}

Consider the maximal solution $u_\Omega$ of the Liouville equation
\begin{equation}
\label{eq:L} -\Delta u + 4 e^{2u}=0
\end{equation}
in $\Omega$. It is known that $u_\Omega$ is the supremum of all
solutions of the Dirichlet problem with smooth boundary data, and
that it is equivalent to $-\ln(2d)$ near the boundary.\footnote{For background
results and applications, see \emph{e.g.}
\cite{BF,K,LN,SK-be,MV,O}. We merely recall that if $\Omega$ is
simply connected, $e^{-u}$ is the ``mapping radius,'' or
``conformal radius'' function of $\Omega$.} Even though the
equation is formally the Euler-Lagrange equation derived from the
Lagrangian $L[u] :=|\nabla u|^2+4e^{2u}$, a direct variational
approach is impossible, because $L[u_\Omega]\not\in L^1(\Omega)$.
Nevertheless, Fuchsian Reduction \cite{SK-note,SK-be,SK-gaeta} enables one
to decompose $u_\Omega$ into an explicit singular part and a $C^1$
function:
\[u_\Omega=v+w,
\]
where the following properties hold
\begin{enumerate}
  \item[(P1)] $w\in C^{1+\alpha}(\overline\Omega)\cap C^2(\Omega)$;
  \item[(P2)] $w=O(d)$ as $d\to 0$;
  \item[(P3)] $e^v=O(1/d)$ as $d\to 0$;
  \item[(P4)] $r[v]:= -\Delta v + 4 e^{2v}=O(1/d)$ as $d\to 0$.
\end{enumerate}
One may, for instance, take $v=-\ln(2d)$ \cite{SK-be}.

For $\phi\in H^1_0$, let us define
\begin{equation}  R[\phi,v] :=
\int_\Omega |\nabla\phi|^2+4e^{2v}[e^{2\phi}-1-2\phi]+2r[v]\phi,
\end{equation}
which is well-defined thanks to eq.~(\ref{eq:exp}), properties
(P3-P4), and Hardy's inequality. We then have a variational
characterization of $u_\Omega$:
\begin{theorem}\label{th:vc-2d}
The infimum
\[ \mathop{\rm Inf} \{ R[\psi-v,v] : \psi\in v+ H^1_0(\Omega)\}
\]
is attained precisely for $\psi=u_\Omega$.
\end{theorem}
Since $v$ is given, this provides a characterization of $u_\Omega$.
\begin{remark}
The result may be stated equivalently as follows: if $\phi\in
H^1_0$,
\begin{equation}\label{eq:equiv}
R[\phi,v]\geq R[w,v],
\end{equation}
with equality if and only if $\phi=w$.
\end{remark}
As a consequence of the variational characterization, we derive a
new global \emph{a priori} bound on $u_\Omega$:
\begin{corollary}\label{cor:4}
The maximal solution $u_\Omega$ of Liouville's equation satisfies
\[\|u_\Omega+\ln(2d)\|_{H^1_0}\leq 2H\|\Delta d\|_{L^2(\Omega)}.
\]
\end{corollary}

\subsection{Organization of the paper}

Section \ref{sec:2} states and proves two results (Theorems
\ref{th:1} and \ref{th:qp}), of which Theorems \ref{th:1-2d} and
\ref{lem:1-2d} are special cases. Th.~\ref{th:vc-2d} and
Cor.~\ref{cor:4} are both proved in Section \ref{sec:v}. The
proofs of Section \ref{sec:2} require the construction of a
partition of unity with special properties, carried out in section
\ref{sec:w}. Comments on the rationale leading to the renormalized
functional are given as concluding remarks.


\section{Hardy-Trudinger inequalities}
\label{sec:2}

In this Section, $\Omega\subset\RR^N$ is an arbitrary domain with
$\pa\Omega\neq\emptyset$. The $N$-dimensional analogues of
Theorems \ref{th:1-2d} and \ref{lem:1-2d} are stated in \S\
\ref{ssec:2-1}, and proved in \S\S\ \ref{ssec:2-2} and
\ref{ssec:2-3} respectively.

\subsection{A synthesis of Hardy's and Trudinger's inequalities}
\label{ssec:2-1}

Let $N'=N/(N-1)$ and define
\[ \Phi_N(u):= \sum_{k\geq N-1} \frac{|u|^{kN'}}{k!},
\]
and, for $1\leq p<\infty$,
\[ M_p(u) :=\left(\int_\Omega
|\nabla u|^p+\frac{|u|^p}{\delta^p} \right)^{1/p}.
\]
We prove:
\begin{theorem}\label{th:1}
If $\pa\Omega\neq\emptyset$, there are constants $c_1$ and $c_2$,
which only depend on the dimension $N$, such that, for any $u\in
W^{1,N}_0(\Omega)$ with $M_N(u)=1$,
\[ \int_{\Omega}
\frac{\Phi_N(u/c_1)}{\delta^{N}}\leq c_2.
\]
\end{theorem}
\begin{remark}
Note that no smoothness or boundedness assumptions on $\Omega$ are
required, and that $M_N(u)$ is not necessarily equivalent to the
$W^{1,N}_0$ norm. This result implies Th.~\ref{lem:1-2d}.
\end{remark}
\begin{remark}
If $\Omega$ is bounded and Lipschitz, Hardy's inequality holds,
and we claim that $\Phi_N(u)$ is integrable for any $u\in
W^{1,N}_0(\Omega)$: write $u=f+g$ where $f$ is smooth with compact
support; since $(|f|+|g|)^{kN'}\leq 2^{kN'}(|f|^{kN'}+|g|^{kN'})$,
we have $\Phi_N(f+g)\leq\Phi_N(2f)+\Phi_N(2g)$. The result follows
if $g$ is small in $W^{1,N}_0$. For $N=2$ and $u\in
H^1_0(\Omega)$, we recover Th.~\ref{th:1-2d}.
\end{remark}

\subsection{An auxiliary result}
\label{ssec:2-2}

Let $1\leq p \leq N$, $p^*=Np/(N-p)$ if $p<N$ (resp.~$p^*=+\infty$
if $p=N$). We prove:
\begin{theorem}\label{th:qp}
If $N\geq 2$, $\Omega$ is a domain in $\RR^N$, $1\leq q<\infty$
and $1\leq p<q\leq p^*$, there is a constant $\Sigma_q(N,p)$ such that
\[ \left(\int_\Omega\frac{|u|^q}{\delta^N}\right)^{1/q}\leq
\Sigma_q\left(\int_\Omega\frac{|\nabla
u|^p}{\delta^{N-p}}+\frac{|u|^p}{\delta^N}\right)^{1/p}
\]
for any $u\in W^{1,p}(\Omega)$.
\end{theorem}
\begin{remark}
In general, the right-hand side may be infinite. If $p=N$ and
Hardy's inequality holds in $\Omega$, the r.h.s.~is finite for
$u\in W_0^{1,p}(\Omega)$.
\end{remark}
\subsection{Proof of Theorem \ref{th:qp}}
\label{ssec:2-3}

{\bf Step 1: Partition of unity.} Denote by $Q(x,s)$ the cube of
center $x$ and side $s$. We prove in section \ref{sec:w} that
there is a smooth partition of unity $(\phi_k)_{k\geq 0}$ in
$\Omega$ with the following properties: \label{ssec:PU}
\begin{enumerate}
\item[(PU1)]
For every $k$, $\phi_k$ is supported in a cube
$Q_k=Q(x_k,s_k)\subset\Omega$.
\item[(PU2)]
For every $k$, $0\leq\phi_k\leq 1$ and $|\nabla\phi_k|\leq
c_3/s_k$, where $c_3$ only depends on $N$.
\item[(PU3)]
There are two positive constants $\lambda$ and $\mu$ such that, on
$\supp\phi_k$, $\lambda \leq \delta/s_k\leq\mu$.
\item[(PU4)]
There is a number $P$ which only depends on the dimension $N$,
such that, for every $x\in\Omega$, at most $P$ among the numbers
$\phi_k(x)$ are non-zero.
\end{enumerate}
Simple consequences of these properties are:
\begin{enumerate}
\item
For every $q\geq 1$, $\sum_k\phi_k^q\leq 1$.
\item $\sum_k \chi_{Q_k}\leq P(N)$,
where $\chi_{Q_k}$ denotes the characteristic function of $Q_k$.
\item $(\sum_k\phi_k)^q\leq P^q\sum_k\phi_k^q$. Indeed, for any
$x$, $\phi_k(x)\neq 0$ for at most $P$ values of $k$, so that
$(\sum_k\phi_k(x))^q\leq P^q\max_k\phi_k(x)^q$.
\end{enumerate}
We also need an elementary observation: for any collection of
nonnegative numbers $b_k$, and any $r\geq 1$,
\begin{equation}\label{eq:r}
\sum_k b_k^r \leq (\sum_k b_k)^r.
\end{equation}
This may be seen for finite sequences by induction, starting from
the inequality: $x^r+y^r\leq (x+y)^r$. Recall also that
$(x+y)^r\leq 2^r(x^r+y^r)$.

{\bf Step 2: Decomposition of $u$.} For any $u\in
W^{1,p}(\Omega)$, we have
\begin{align*}
\int_\Omega \frac{|u|^q}{\delta^{N}} &=
    \int_\Omega |\sum_k u\phi_k|^q \delta^{-N} \\
&\leq P^q\sum_k \int_{Q_k} |u\phi_k|^q \delta^{-N}\\ &\leq
P^q\sum_k (\lambda s_k)^{-N}\|u\phi_k\|_{L^q(Q_k)}^q.
\end{align*}
Write $Q(s)$ for $Q(0,s)$, and let $S_q=S_q(N,p)$ denote the norm
of the embedding of $W^{1,p}_0(Q(1))$ into $L^q(Q(1))$. If $v\in
W^{1,p}_0(Q(s))$, one finds, applying the Sobolev inequality to
$v(s x)\in W^{1,p}_0(Q(1))$,
\begin{equation}
\|v\|_{L^q(Q(s))} \leq S_qs^{\frac Nq+1-\frac Np}\|\nabla
v\|_{L^p(Q(s))}.
\end{equation}
It follows that, for every $k$,
\[ (\lambda s_k)^{-N}\|u\phi_k\|_{L^q(Q_k)}^q
\leq S_q^q\lambda^{-N} s_k^{q(p-N)/p}\|\nabla(u\phi_k)\|_{L^p(Q_k)}^q.
\]
It follows that
\begin{align*}
\left(\int_\Omega \frac{|u|^q}{\delta^{N}}\right)^{p/q} & \leq
(PS_q)^p \left(\sum_k \lambda^{-N} s_k^{q(p-N)/p}
\|\nabla(u\phi_k)\|_{L^p(Q_k)}^q\right)^{p/q}\\ & \leq
(PS_q\lambda^{-N/q})^p \sum_k
s_k^{p-N}\|\nabla(u\phi_k)\|_{L^p(Q_k)}^p\\ & \leq
(2PS_q\lambda^{-N/q})^p \sum_k  s_k^{p-N}\left[\|\phi_k\nabla
u\|_{L^p(Q_k)}^p+\|u\nabla\phi_k\|_{L^p(Q_k)}^p\right],
\end{align*}
where we used eq.~(\ref{eq:r}) to obtain the second inequality.
Now,
\[\int_\Omega\sum_k s_k^{p-N}\phi_k^p|\nabla u|^p\leq
\mu^{N-p}\int_\Omega\frac{|\nabla u|^p}{\delta^{N-p}},
\]
and
\begin{eqnarray*}
\sum_k s_k^{p-N}\|u\nabla\phi_k\|_{L^p(Q_k)}^p & \leq &
c_3^p\int_\Omega(\sum_k \chi_{Q_k}(x))\frac{|u|^p}{s_k^N}\\ & \leq
& Pc_3^p\mu^N\int_\Omega \frac{|u|^p}{\delta^N}.
\end{eqnarray*}
We have therefore the desired inequality, with
\begin{equation}\label{eq:ct}
\Sigma_q=2PS_q(N,p)\lambda^{-N/q}[\mu^{N-p}+Pc_3^p\mu^N]^{1/p}.
\end{equation}
This completes the proof.

\subsection{Proof of Theorem \ref{th:1}}
\label{ssec:3}

We now consider the case $p=N$, so that $p^*=+\infty$, and $q$ can
take arbitrarily large values. From \cite[lemma 7.12 and
eq.~(7.37)]{GT}, it follows that
\[ S_q(N,N)\leq (\omega_Nq)^{1-1/N+1/q}\text{ if } q\geq N.
\]
If $q\geq N-1$, we have $N'q\geq N$, and therefore,
\[ S_{N'q}(N,N)^{N'q}\leq (N'q\omega_N)^{q+1}\text{ if } q\geq N-1.
\]
For any $c_1>0$, we therefore find
\[
\int_\Omega\sum_{q\geq N-1} \frac{|u|^{N'q}}{q!c_1^{N'q} \delta^N}
\leq c_2:=\sum_{q\geq N-1} \lambda^{-N}N'\omega_N
(\frac{N'\omega_NA}{c_1^{N'}})^q\frac{q^q}{(q-1)!},
\]
where $A=\{2P[1+P(c_3\mu)^N]^{1/N}\}^{N'}$. The series defining $c_2$
converges if $c_1^{N'}>e\omega_N N'A$.

This completes the proof.


\section{Variational characterization of solutions with boundary blow-up}
\label{sec:v}

\subsection{Proof of theorem \ref{th:vc-2d}}

Let $\phi\in H^1_0(\Omega)$. We wish to prove inequality
(\ref{eq:equiv}). First observe that since $u_\Omega=v+w$ solves
Liouville's equation,
\begin{equation}\label{eq:bb-1}
\Delta w=-\Delta v + 4e^{2v}+4e^{2v}(e^{2w}-1)
           = r[v]+4e^{2v}(e^{2w}-1).
\end{equation}
It follows from (P2-P4) that
\[ \Delta w=O(1/d).
\]
Since $w\in C^1_0(\overline\Omega)\cap C^2(\Omega)$, we have,
for any $\psi\in C^1_0(\overline\Omega)$ and any $\ep>0$ small
enough,
\[ \int_{d>\ep} \nabla\psi\cdot\nabla w=
\int_{d=\ep}\psi(\nabla w\cdot\nabla d)\;ds -
 \int_{d>\ep} \psi\Delta w.
 \]
Letting first $\ep\to 0$, and then approximating $\phi$ by $\psi$
in the $H^1_0$ norm, we find
\begin{equation}\label{eq:bb-2}
\int_\Omega\nabla\phi\cdot\nabla w+\int_\Omega \phi\Delta w= 0.
\end{equation}
Using equation (\ref{eq:bb-1}), we find
\begin{equation}\label{eq:bb-3}
\int_\Omega\nabla\phi\cdot\nabla w +4e^{2v}(e^{2w}-1)\phi+r[v]\phi
=0.
\end{equation}
Since
\[R[\phi+w,v]-R[w,v]=\int_\Omega|\nabla\phi|^2+2\nabla\phi\cdot\nabla
w+4e^{2v}[e^{2w+2\phi}-e^{2w}-2\phi]+2r[v]\phi,
\]
we find
\begin{equation}\label{eq:bb-4}
R[\phi+w,v]-R[w,v]=
 \int_\Omega|\nabla\phi|^2+4e^{2(v+w)}(e^{2\phi}-1-2\phi),
\end{equation}
which is manifestly nonnegative, and vanishes precisely if
$\phi=0$. Q.E.D.
\begin{remark}
Since $v+w=u_\Omega$ and $r[u_\Omega]= 0$, the r.h.s.~of equation
(\ref{eq:bb-4}) is equal to $R[\phi,u_\Omega]$.
\end{remark}

\subsection{Proof of Corollary \ref{cor:4}}

Property (P4) and Hardy's inequality ensure that there is a
constant $K[v]$ such that
\[-\int_\Omega 2r[v]\phi \leq K[v] \|\phi\|_{H^1_0}.
\]
Expressing that $R[w,v]\leq R[0,v]=0$, we find
\[\|w\|_{H^1_0}\leq K[v].
\]
If $v=-\ln(2d)$, one finds $r[-\ln(2d)]=(\Delta d)/d$. H\"older's
and Hardy's inequalities yield $K[v]\leq 2H \|\Delta d\|_{L^2}$.
Since $w=u_\Omega-v$, the announced \emph{a priori} $H^1$ bound on
$u_\Omega$ follows.


\section{Construction of the partition of unity}
\label{sec:w}

We construct the partition of unity used in the proof of Theorem
\ref{th:qp}; the basic ideas go back to Whitney \cite{W}, and many
variants may be found in the literature, see e.g.~\cite{S,Mz}.

First of all, choose two constants $\eta$ and $\eta'$ such that
\[\eta/\sqrt N>\eta'>1.
\]
Recall that $Q(x,s)$ is the closed cube of center $x$ and side $s$. For
any $\sigma>0$, we write $Q_\sigma(x,s)$ for $Q(x,\sigma s)$. The
following observation will be useful: for any $x\in\Omega$,
\[ \delta(x)>\frac s2 \sqrt N \Rightarrow Q(x,s)\subset\Omega
                \Rightarrow \delta(x)>\frac s2.
\]
Conversely,
\[ \delta(x)\leq \frac s2  \Rightarrow Q(x,s)\not\subset\Omega
                \Rightarrow \delta(x)\leq\frac s2 \sqrt N.
\]

\subsection{Covering by dyadic cubes}

For $k\in\mathbb Z$, let ${\mathcal F}_k$ denote the family of
closed cubes of the form
\[[0,2^{-k}]^N+(m_1,\dots,m_N)2^{-k},
\]
where the $m_j$ are signed integers. The cubes in ${\mathcal
F}_{k+1}$ are obtained by dyadic subdivision of the cubes in
${\mathcal F}_k$; in particular, for any $k$, every
$Q(x,s)\in\mathcal{F}_{k+1}$ is included in a unique cube $\tilde
Q(\tilde x,\tilde s)\in\mathcal{F}_{k}$, with $\tilde s=2s$. Since
$Q$ is obtained by dyadic division of $\tilde Q$, $\tilde x$ must
be a vertex of $Q$; it follows that $|x-\tilde x|=\frac12 s\sqrt
N$.

Let ${\mathcal F}=\bigcup_{k=-\infty}^{+\infty}{\mathcal F}_k$.
Define a set $\mathcal Q\subset\mathcal F$ as follows:
$Q\in\mathcal Q$ if and only if
\begin{equation}
Q_\eta\subset\Omega \text{ and } \tilde
Q_\eta\not\subset\Omega.
\end{equation}
$\mathcal F$ is not empty since $\pa\Omega\neq\emptyset$ by
assumption.
\begin{lemma}
$\Omega$ is the union of the cubes $Q\in\mathcal Q$.
\end{lemma}
\begin{proof}
Let $y\in\Omega$. Consider the set of numbers $k$ for which there
is a cube $Q(x,2^{-k})\in\mathcal{F}_k$ which contains $y$ and
which satisfies $Q_{\eta}\subset\Omega$. This set is not empty: if
$k$ is large enough, we have $\delta(x)\geq\delta(y)-\frac12
2^{-k}\sqrt N>\frac12 2^{-k}\eta\sqrt N$, and
$Q_{\eta}\subset\Omega$. It is bounded below because
$\pa\Omega\neq\emptyset$. Let $k_0$ be the smallest integer in
that set, and consider a cube $Q$ with the above property with
$k=k_0$. Since $k_0$ is minimal, $\tilde Q_\eta\not\subset\Omega$.
Therefore, $y\in Q\in\mathcal Q$, as desired.
\end{proof}
Since, for each $Q\in\mathcal Q$, $Q\subset
Q_{\eta'}\subset\Omega$, we have \emph{a fortiori}
\[\Omega=\bigcup_{Q\in\mathcal Q} Q_{\eta'}.
\]

\subsection{Properties of the cube decomposition}

We now prove that the covering of $\Omega$ by the cubes
$Q_{\eta'}$ has the additional property that, on each of them, the
function $\delta$ is comparable to the side of $Q$:
\begin{lemma}\label{lem:2}
There are positive constants $c_4$ and $c_5$, independent of
$\Omega$, such that, if $Q(x,s)\in\mathcal Q$ and $y\in
Q_{\eta'}$, then
\[c_4\leq \frac{\delta(y)}{s} \leq c_5.
\]
\end{lemma}
\begin{proof}
Since $Q_\eta\subset\Omega$, $\delta(x)>\eta s/2$. If $y\in
Q_{\eta'}$, $\delta(y)\geq\delta(x)-\frac12 \eta's\sqrt
N>\frac12(\eta-\eta'\sqrt N)s$. Therefore,
\[\frac12(\eta-\sqrt N{\eta'})<\frac{\delta(y)}{s}.
\]
To establish an upper bound, we first estimate $\delta(x)/s$.
Since $Q(x,\eta s)\subset\Omega$ and $\tilde Q(\tilde x,2\eta
s)\not\subset\Omega$, $\delta(x)>\frac12 \eta s$ and
$\delta(\tilde x)\leq\frac12 (2\eta s) \sqrt N$. Therefore,
\[\delta(x)\leq \delta(\tilde x)+|x-\tilde x|\leq (\eta+\frac12)s\sqrt N.
\]
Therefore, if $Q(x,s)\in\mathcal{Q}$,
\begin{equation}\label{eq:w1}
\frac12 \eta < \frac{\delta(x)}s \leq (\eta+\frac12)\sqrt N.
\end{equation}
We now estimate $\delta(y)$:
\[
\delta(y)\leq\delta(x)+\frac12 \eta's\sqrt N\leq
(\eta+\frac12+\frac12\eta')s\sqrt N.
\]
It follows that, $\delta(y)/s$ lies between positive bounds which
do not depend on $\Omega$, as desired.
\end{proof}

\subsection{Finite intersection property}

The cubes $Q_{\eta'}$ are not disjoint; nevertheless, two such
cubes may intersect only if the ratio of their sides lies between
fixed bounds; this implies that a given cube may only intersect
finitely many others. More precisely,
\begin{lemma}\label{lem:fip}
There is a constant $c_6$, independent of $\Omega$, such that, if
$Q(x,s)$ and $Q'(x',s')$ belong to $\mathcal Q$, and if
$Q(x,\eta's)\cap Q'(x',\eta's')\neq\emptyset$, with $s<s'$,
necessarily
\[ s < s' < c_6 s.
\]
Furthermore, there is a number $P$, independent of $\Omega$, such
that given $Q\in\mathcal Q$, there are at most $P$ cubes
$Q'\in\mathcal Q$ such that $Q_{\eta'}\cap
Q'_{\eta'}\neq\emptyset$.
\end{lemma}
\begin{proof}
Let $y\in Q_{\eta'}\cap Q'_{\eta'}$. We have $|x-x'|\leq
|x-y|+|y-x'|\leq \frac12 (s+s')\eta'\sqrt N$. Since $Q$ and $Q'$
both belong to $\mathcal Q$, inequality (\ref{eq:w1})
yields $\delta(x')>\frac12 \eta s'$ and $\delta(x)\leq
(\frac12+\eta)s\sqrt N$. Since $\delta$ is 1-Lipschitz,
\[\frac12 \eta s'<\delta(x')\leq \delta(x)+|x-x'|\leq
[(\frac12+\eta+\frac12\eta')s+\frac12\eta's']\sqrt N.
\]
It follows that $s'<c_6s$ with
\[ c_6=\frac{2\eta+1+\eta'}{\eta-\eta'\sqrt N}\sqrt N.
\]
Therefore, if $s=2^{-k}$ and $s'=2^{-k'}$, we have $|k-k'|\leq
J_1:=\ln c_6/\ln 2$. On the other hand, $|x-x'|\leq\frac12
(s+s')\eta'\sqrt N\leq J_2s$, with $J_2=\frac12 (1+c_6)\eta'\sqrt N$.
Scaling the cubes by the factor $1/s$, and taking $x$ as the
origin of coordinates, we are led to the question: ``how many
cubes $Q(x',2^j)$ can one find subject to the restrictions:
$|j|\leq J_1$ and $|x'|\leq J_2$?'' The answer is a finite number
$P$, which depends only on $N$, $\eta$ and $\eta'$.
\end{proof}

\subsection{Constructing the partition of unity}

Take a smooth function $\varphi(x)$ with support in $Q(0,\eta')$,
equal to 1 on $Q(0,1)$, and such that $0\leq\varphi(x)\leq 1$ for
all $x$. Consider
\[\psi(x)=\sum_{Q=Q(x_Q,s_Q)\in\mathcal Q}\varphi(\frac{x-x_Q}{s_Q}).
\]
Since the cubes $Q\in\mathcal Q$ cover $\Omega$, $\psi(x)\geq 1$
for all $x\in\Omega$. Since any point belongs to at most $P$ of
the cubes $Q_{\eta'}$, we have $\psi(x)\leq P$ for all
$x\in\Omega$. It follows that the functions
$\phi_Q(x):=\varphi((x-x_Q)/{s_Q})/\psi(x)$ form a partition of
unity, and are supported in the cubes $Q_{\eta'}$.

We now prove that properties (PU1--4) of section \ref{ssec:PU}
hold. (PU1) and (PU2) are immediate. If $x\in Q_{\eta'}$, Lemma
\ref{lem:2} ensures that $\delta(x)/s_Q$ is bounded above and
below by positive bounds which only depend on $N$. This proves
property (PU3). Finally, since $\varphi((x-x_Q)/{s_Q})$ is
supported in $Q_\eta'$, we find, thanks to lemma \ref{lem:fip},
that (PU4) holds.

This completes the construction of the partition of unity.

\section{Concluding remarks}

The boundary blow-up problem appears at first sight to be beyond
the reach of critical-point theory. We have shown that it is in
fact equivalent to the minimization of a convex functional in
$H^1_0$.

How is this functional related to the usual functional
$E[u]:=\int_\Omega L[u]$? Take $v=-\ln(2d)$ to fix ideas. Even
though $E[u]$ is infinite for $u=u_\Omega$, it is easy to see
that, if $\phi$ is smooth and sufficiently flat near the boundary,
$\tilde E[\phi]:=\int_\Omega L[\phi+v]-L[v]$ is well-defined. But
this does not provide a satisfactory variational principle for two
reasons: (i) $u_\Omega-v$ is not very flat at the boundary: it is
only $O(d)$; (ii) $\tilde E[\phi]$ is not well-defined if $\phi$
is merely $O(d)$ (indeed, the term $2\nabla\phi\cdot\nabla v$ is
not necessarily integrable). However, subtracting $\int_\Omega
2\div(\phi\nabla v)$ from $\tilde E[\phi]$, and using the equation
satisfied by $v$, one recovers an expression equivalent to our
functional $R$. Note that the integral of this divergence term is
not zero, even if $\phi$ is smooth, because $\nabla v$ blows up at
the boundary.

Similar ideas are applicable in other situations, as soon as one
can decompose the singular solution into the sum $v+w$ of an
explicit singular function, and a remainder of controlled
regularity. Fuchsian Reduction provides a systematic procedure for
achieving this decomposition.

\end{document}